\documentclass{birkjour}
\usepackage{amsmath}
\usepackage{amssymb}
\usepackage{url}

\usepackage{yhmath} 

\usepackage{srctex} 

\usepackage{color}

\newcommand{\bbH}{\ensuremath{\mathbb{H}}}

\newcommand{\ii}{\mathrm{i}} 
\newcommand{\jj}{\mathrm{j}}
\newcommand{\kk}{\mathrm{k}}

\newcommand{\reverse}[1]{\widetilde{#1}}
\newcommand{\gradeinverse}[1]{\wideparen{#1}}
\newcommand{\cliffordconjugate}[1]{\widetilde{\wideparen{#1}}}

\newcommand{\cl}[2]{\ensuremath{\mathit{Cl}_{#1,#2}}} 

\newcommand{\e}[1]{\mathbf{e}_{#1}} 
\newcommand{\ba}{\mathbf{a}} 

\begin{document}


\title[A note on solution of $ax+xb=c$ by Clifford algebras]
 {A note on solution of $ax+xb=c$ by Clifford algebras}

\author{A.~Dargys}
\address{%
Center for Physical Sciences and Technology,\br Semiconductor
Physics Institute,\br Saul{\.e}tekio 3, LT-10257 Vilnius,
Lithuania} \email{adolfas.dargys@ftmc.lt}

\author{A.~Acus${}^\dagger$}
\address{Institute of Theoretical Physics and Astronomy,\br
Vilnius University,\br Saul{\.e}tekio 3, LT-10257 Vilnius,
Lithuania} \email{arturas.acus@tfai.vu.lt} \thanks{$\dagger$ Corresponding author}



\date{October, 2018}

\begin{abstract}
The coordinate-free solutions of the multivector equation
$ax+xb=c$ are discussed and presented for the Clifford algebras
\cl{p}{q} when $p+q\le 3$.
\end{abstract}

\maketitle

\section{Introduction\label{sec:1}}

The noncommutative linear equation in  quaternionic variables or
in a matrix form,
\begin{equation}\label{sylv}
ax+xb=c,\quad a,b,c,x\in\bbH,
\end{equation}
where  $x$ is the unknown, plays a central role in many areas of
applied mathematics, in particular, in signal processing and
control theory. There is a vast literature on properties and
methods  of solutions of~\eqref{sylv}. The articles
\cite{Niven1941,Eilenberg1944,Turner2006,Janovska2008,Helmstetter2012,schwartz2013},
for example, may be considered as representative.

The quaternion  $a=a_0+a_1\ii+a_2\jj+a_3\kk=a_0+\ba$ consists of
scalar $a_0$ and vector $\ba$  parts. The latter is characterized
by imaginary elements that satisfy $\ii^2=\jj^2=\kk^2=-1$, and
$\ii\jj\kk=-1$. The conjugate quaternion,
$\bar{a}=a_0-a_1\ii-a_2\jj-a_3\kk=a_0-\ba$,  allows to find the
inverse $a^{-1}=\overline{a}/(a\overline{a})$. The sum $a+\bar{a}$
and the product $a\overline{a}$  are the scalars, commutativity of
which play the central role in quaternionic problems, including
the solution of~\eqref{sylv} in a coordinate-free way.

The quaternion algebra is isomorphic to Clifford \cl{0}{2} algebra
with $\ii\simeq\e{1}$, $\jj\simeq\e{2}$ and $\kk\simeq\e{12}$,
where $\e{1}$ and $\e{2}$ are basis vectors and $\e{12}$ is the
bivector, $\e{12}\equiv\e{1}\e{2}$. The quaternionic conjugation
in \cl{0}{2} is replaced by Clifford conjugation
$\cliffordconjugate{a}$ which is defined as a simultaneous action
of the grade inversion $\gradeinverse{a}$ (which maps $\e{i}\to
-\e{i}$) and the reversion $\reverse{a}$ (which corresponds to
$\e{12}\to\e{21}=-\e{12}$). In a search of possible solutions of
Sylvester equation~\eqref{sylv} in the Clifford algebra  it is
convenient to replace the Clifford conjugation by grade negation
operation as $\cliffordconjugate{a}=a_{\bar{1},\bar{2}}$, where
the bar over subscripts designates which of the grades have
opposite signs. Thus, $\bar{1}$ changes signs of all vectors to
opposite and similarly $\bar{2}$ changes signs of all
bivectors~\cite{Acus2018}. In Sect.~\ref{sec:2} we shall remind
how the Sylvester equation~\eqref{sylv} can be solved in a
coordinate-free way and in Sect.~\ref{sec:3} the method is applied
to solve~\eqref{sylv} for the Clifford algebras of dimension
$2-4$.

\section{Coordinate-free solution of quaternionic Sylvester equation\label{sec:2}}
At first let us remind how the equation~\eqref{sylv} is solved in
a coordinate-free way, where for a moment the overbar means the
quaternionic conjugate.

\noindent1) Multiply the equation~\eqref{sylv} from left by $a$ to
get
\begin{equation}
a^2x+axb=ac.
\end{equation}
2) Then multiply  from right by $\bar{b}$ to get
\begin{equation}
ax\bar{b}+xb\bar{b}=c\bar{b}.
\end{equation}
3) Since both  equations are linear with respect to unknown  they
can be added,
\begin{equation}
(a^2x+xb\bar{b})+(axb+ax\bar{b})=ac+c\bar{b}.
\end{equation}
4) Since $b\bar{b}$ and $b+\bar{b}$ are scalars, we can write
\begin{equation}
(a^2x+b\bar{b})x+a(b+\bar{b})x=ac+c\bar{b},
\end{equation}
from which we have the solution:
\begin{equation}\label{solSylvA}
x=\big[(a^2+b\bar{b})+a(b+\bar{b})\big]^{-1}(ac+c\bar{b}).
\end{equation}
This formula coincides with the quaternionic expression given
in~\cite{Janovska2008} on p.225.

If instead, the Sylvester equation~\eqref{sylv} is at first
multiplied from the right by $b$ and from the left by $\bar{a}$ then one
gets a mirror expression of~\eqref{solSylvA},
\begin{equation}\label{solSylvB}
x=(\bar{a}c+cb)[(b^2+\bar{a}a)+b(a+\bar{a})]^{-1}.
\end{equation}
From both equations, \eqref{solSylvA} and \eqref{solSylvB}, one
can see that the solution of Sylvester equation has been reduced
to a problem of quaternionic inverse, or as we shall show in the
next section to the inverse of a multivector in Clifford algebras.

\section{Solution of $ax+xb=c$ by Clifford algebra\label{sec:3}}

The solution of quaternionic Sylvester equation in coordinate-free
form has been easy because it was possible  to form two scalars,
$b+\bar{b}$ and $b\bar{b}$, that commute with all quaternions.
More generally, a center of  algebra $\mathcal{A}$ may be
introduced that consists of those elements $s\in\mathcal{A}$ which
commute with all elements of the considered algebra,
\begin{equation}
  \mathrm{Cen}(\mathcal{A})=\{s\in\mathcal{A}\;|\;as=sa, \forall
a\in\mathcal{A}\}.
\end{equation}
For example, we can solve the Sylvester equation in a similar
manner for \cl{3}{0} algebra if a more general center consisting
of scalar and pseudoscalar is constructed from the multivector
$a$, i.~e., $\mathrm{Cen}_1=a+\bar{a}\in \mathrm{Cen}(\cl{3}{0})$
and $\mathrm{Cen}_2=a\bar{a}\in \mathrm{Cen}(\cl{3}{0})$.

More generally, one may try  to construct solutions in the
form~\eqref{solSylvA} and \eqref{solSylvB} if instead of Clifford
conjugation one introduces new involutions, or operations (in the
following such operations will be denoted by overbar too) that
simultaneously satisfy either
$(b+\bar{b})=\mathrm{Cen}_1(\cl{p}{q})$ and
$(b\bar{b})=\mathrm{Cen}_2(\cl{p}{q})$, or in case of multivector
$a$, correspondingly, $(a+\bar{a})=\mathrm{Cen}_1(\cl{p}{q})$ and
$(a\bar{a})=\mathrm{Cen}_2(\cl{p}{q})$. Then  repeating the
procedure described  in Sect.~\ref{sec:2} one can automatically
get the solutions \eqref{solSylvA} and \eqref{solSylvB}. 

Table~\ref{centergrades} shows those Clifford algebras (CA) where
the CA Sylvester equation
\begin{equation}\label{sylvGA}
ax+xb=c,\quad a,b,c,x\in\cl{p}{q},
\end{equation}
is solvable in the form \eqref{solSylvA} or~\eqref{solSylvB}. The
solutions in Clifford algebras (contrary to the quaternion case
where the inverse exists always if $a\ne 0$) have singularity when
$(a^2+b\bar{b})+a(b+\bar{b})=0$, or
$(b^2+\bar{a}a)+b(a+\bar{a})=0$.  If $b=a$ these equations reduce
to $x=\big[2(a+\bar{a})\big]^{-1}(c+a^{-1}c\bar{a})$, which are
singular when  $a+\bar{a}=0$, i.~e, when both the scalar and the
pseudoscalar are equal to zero).

\begin{table}
\begin{center}
\caption{Clifford algebras (CA) that allow coordinate-free
solution of equation~\eqref{sylvGA} in the form~\eqref{solSylvA}
or ~\eqref{solSylvB}.  The last two columns represent CA centers.
$I_3$ is the corresponding  pseudoscalar.}
\begin{tabular}{c|ccc}\label{centergrades}
CA &$\bar{a}$ &$\mathrm{Cen}_1=a+\bar{a}$ &$\mathrm{Cen}_2=a\bar{a}$ \\
\hline
\cl{2}{0} &$a_{\bar{1},\bar{2}}$&$2a_0$ &$a_0^2-a_1^2-a_2^2+a_3^2$ \\
\cl{1}{1} &$a_{\bar{1},\bar{2}}$&$2a_0$ &$a_0^2-a_1^2+a_2^2-a_3^2$ \\
\cl{0}{2} &$a_{\bar{1},\bar{2}}$&$2a_0$ &$a_0^2+a_1^2+a_2^2+a_3^2$
\\[5pt]
\cl{3}{0} &$a_{\bar{1},\bar{2}}$&$2(a_0+a_{123}I_3)$ &$a_0^2-a_1^2-a_2^2-a_3^2+a_{12}^2+a_{13}^2+a_{23}^2-a_{123}^2$ \\
 & & &$+2(a_0a_{123}-a_1a_{23}+a_{2}a_{13}-a_{3}a_{12})I_3$ \\
 \cl{2}{1} &$a_{\bar{1},\bar{2}}$&$2(a_0+a_{123}I_3)$ &$a_0^2-a_1^2-a_2^2+a_3^2+a_{12}^2-a_{13}^2-a_{23}^2+a_{123}^2$ \\
 & & &$+2(a_0a_{123}-a_1a_{23}+a_{2}a_{13}-a_{3}a_{12})I_3$ \\
 \cl{1}{2} &$a_{\bar{1},\bar{2}}$&$2(a_0+a_{123}I_3)$ &$a_0^2-a_1^2+a_2^2+a_3^2-a_{12}^2-a_{13}^2+a_{23}^2-a_{123}^2$ \\
 & & &$+2(a_0a_{123}-a_1a_{23}+a_{2}a_{13}-a_{3}a_{12})I_3$ \\
 \cl{0}{3} &$a_{\bar{1},\bar{2}}$&$2(a_0+a_{123}I_3)$ &$a_0^2+a_1^2+a_2^2+a_3^2+a_{12}^2+a_{13}^2+a_{23}^2+a_{123}^2$\\
 & & &$+2(a_0a_{123}-a_1a_{23}+a_{2}a_{13}-a_{3}a_{12})I_3$
\end{tabular}
\end{center}
\end{table}

\vspace{3mm} \textit{Example}. Given $\cl{3}{0}$ algebra
multivectors, for example $a=3+3 \e{1}+2 \e{13}+5 \e{123}$, $b=3+2
\e{2}+3 \e{3}+2 \e{123}$, and $c=5 \e{1}+3 \e{2}+4 \e{13}+\e{23}$,
let us find the solution of Sylvester equation~\eqref{sylvGA} by
applying formula \eqref{solSylvA}. We obtain: 1) $\bar{b}=3-2
\e{2}-3 \e{3}+2 \e{123}$, 2) $(a^2+b\bar{b})+a(b+\bar{b})=-21+36
\e{1}+28 \e{2}+24 \e{13}+42 \e{23}+84 \e{123}$.  3) The inverse of 3D~multivector~$A$ 
 is $A^{-1}=\frac{C(A
C)_{\bar{3}}}{AC(A C)_{\bar{3}}}$, where $C=A_{\bar{1},\bar{2}}$,
which gives
$\bigl[(a^2+b\bar{b})+a(b+\bar{b})\bigr]^{-1}=\frac{1}{2177719}
(-9807+14436 \e{1}+1708 \e{2}+9624 \e{13}+2562 \e{23}-20748
\e{123})$.  4) Then $(ac+c\bar{b})=1+11 \e{1}+43 \e{2}+4 \e{3}-3
\e{12}-12 \e{13}+32 \e{23}+5 \e{123}$, and finally 5) $x=
\frac{1}{2177719} (359677+601305 \e{1}-155957 \e{2}-436078
\e{3}+209677 \e{12}+1076362 \e{13}-489350 \e{23}+27015 \e{123})$.
One can check that the obtained solution indeed satisfies
equation~\eqref{sylvGA}. In a similar way one can obtain the
solution using expression~\eqref{solSylvB}, though the
intermediate steps will give different multivectors. Since even
subalgebras of larger $n=p+q=4$ algebra  are isomorphic to one of
$n=3$ algebra the method can also be applied, although with
limited capabilities, to CA algebras \cl{1}{3} and \cl{3}{1} that
are very important  in the relativistic physics.

The  grade negation operation allows to detect accidental centers,
if such exist, in the higher grade algebras. We have also
investigated the algebras with $n=4,5$ and $6$ applying for this
purpose the grade negation operation. However, we have been able
to construct only a single CA center rather than two centers that
are required by algorithm in Sect.~\ref{sec:2} to solve the
Sylvester equation~\eqref{sylvGA}.

\end{document}